\documentclass[a4paper,11pt]{amsart}
\usepackage[utf8]{inputenc}
\usepackage{amsaddr}
\usepackage{calrsfs}
\usepackage{graphicx,color}
\usepackage{hyperref}

\numberwithin{equation}{section}

\theoremstyle{definition}

\theoremstyle{plain}
\newtheorem{thm}{Theorem}[section]

\theoremstyle{definition}
\newtheorem{rem}{Remark}[section]

\newcommand{\R}{\mathbb{R}}
\newcommand{\E}{\mathbb{E}}

\newcommand{\F}{\mathcal{F}}

\newcommand{\1}{\mathbf{1}}
\renewcommand{\d}{\mathrm{d}}

\begin{document}
\title[Prediction Law of Mixed Gaussian Volterra Processes]
{Prediction Law of Mixed Gaussian Volterra Processes}

\date{\today}

\author[Sottinen]{Tommi Sottinen}
\address{Department of Mathematics and Statistics, University of Vaasa, P.O. Box 700, FIN-65101 Vaasa, FINLAND}
\email{tommi.sottinen@iki.fi}

\author[Viitasaari]{Lauri Viitasaari}
\address{Department of Mathematics,  University of Helsinki, P.O.Box 68, FIN-00014 University of Helsinki, FINLAND} 
\email{lauri.viitasaari@iki.fi}

\begin{abstract}
We study the regular conditional law of mixed Gaussian Volterra processes under the influence of model disturbances. More precisely, we study prediction of Gaussian Volterra processes driven by a Brownian motion in a case where the Brownian motion is not observable, but only a noisy version is observed. As an application, we discuss how our result can be applied to variance reduction in the presence of measurement errors.
\end{abstract}




\maketitle

\section{Introduction}

We study the regular conditional prediction law of mixed Gaussian Volterra processes. More precisely, we study prediction of a given Gaussian Volterra process in a case where one does not observe the process directly, but instead observes a noisy version of it. Such problems arise naturally in cases when there are measurements errors. We apply our results to variance reduction in such a case. 

Surprisingly, regular conditional laws have not been studied extensively in the literature. On related research we can mention \cite{Shokrollahi-Sottinen-2017}, \cite{Sottinen-Viitasaari-2015}, \cite{Sottinen-Viitasaari-2016b} and \cite{Sottinen-Viitasaari-2017b} where fractional Brownian motion and more generally Gaussian Fredholm processes were considered and applied in the stochastic finance setting.

The rest of the paper is organized as follows. In Section \ref{sec:main} we introduce our setting and formulate and prove our main result. Section \ref{sec:application} is devoted to the application to variance reduction in the case of measurement errors.  

\section{Regular Conditional Prediction Laws}
\label{sec:main}
We consider the prediction of a Gaussian Volterra process $X=(X_t)_{t\ge 0}$ defined as
\begin{equation}\label{eq:gvp}
X_t = \int_0^t k(t,s)\, \d W_s,
\end{equation}
where $W=(W_t)_{t\ge0}$ is a Brownian motion and $k\colon\R^2\to\R$ is a deterministic kernel such that 
$$
\int_0^t k(t,s)^2\, \d s < \infty
$$
for all $t\ge 0$ so that the Wiener integral in \eqref{eq:gvp} is well-defined. We denote by $r\colon\R^2\to\R$ the covariance of $X$ and note that by the representation \eqref{eq:gvp} we have
$$
r(t,s) = \int_0^{t\wedge s} k(t,u)k(s,u)\, \d u,
$$
where we have used the common notation $t\wedge s = \min(t,s)$.

Let $\tilde W=(\tilde W_t)_{t\ge 0}$ be another Brownian motion that is independent of $W$. Let $a,b\in\R$ and denote $W^{a,b} = aW + b\tilde W$. We are interested in the regular conditional law of the process $X$ given by the filtration $\mathbb{F}^{a,b}=(\F^{a,b}_u)_{u\ge 0}$ of the mixed Brownian motion $W^{a,b}$. The parameters $a,b\in\R$ are assumed to be known as well as the kernel $k$. 

We denote by $\1_t$ the indicator
$$
\1_t(s) = \left\{\begin{array}{rl}
1, & \mbox{if } s<t, \\
0, & \mbox{otherwise}
\end{array}
\right.
$$


Since, in the Gaussian setting, the conditional expectation is a linear operator of the observed path, one expects that
$$
\hat m_t^{a,b}(u) = \E\big[X_t\,\big|\,\F_u^{a,b}\big] 
=
\int_0^u \hat k^{a,b}(t,s|u)\, \d W^{a,b}_s,
$$
where $\hat k$ is some suitable kernel. It turns out that this is indeed the case. Moreover, $\hat k^{a,b}(t,s|u)$ is independent of $u$. The precise statement is the following result.

\begin{thm}[Prediction Law]\label{thm:cond_law}
Let $u>0$. The regular conditional law of $X|\F^{a,b}_u$ on $[0,\infty)$ is Gaussian with random mean
\begin{equation}\label{eq:cond_mean}
\hat m_t^{a,b}(u) = \frac{a}{a^2+b^2}\int_0^u k(t,s)\, \d W^{a,b}_s
\end{equation}
and deterministic covariance
\begin{eqnarray}\label{eq:cond_var}
\hat r^{a,b}(t,s|u) &=& 
\int_0^{t\wedge s} \left(1-\frac{a^2}{a^2+b^2}\1_u(v)\right)^2 k(t,v)k(s,v)\, \d v \\ \nonumber
& &
+
\frac{a^2b^2}{(a^2+b^2)^2}\int_0^u k(t,v)k(s,v)\, \d v.  
\end{eqnarray}
\end{thm}

\begin{rem}
Note that here we consider the prediction law on $[0,\infty)$, meaning that we are also interested in the prediction of the past and the present, not only the future. Indeed, predicting the past and the present are also issues, since the process $X$ is not directly observable from $W^{a,b}$. This situation is closely related to Gaussian bridges as studied, e.g., in \cite{Gasbarra-Sottinen-Valkeila-2007} and \cite{Sottinen-Yazigi-2014}.  
\end{rem}

\begin{rem}[Correlation Paramerization]
The prediction with the parameters $(a,0)$ is the same as with the parameters $(1,0)$.  Indeed, in the case $(a,0)$ we do not have the ``disturbing'' Brownian motion $\tilde W$, and the prediction is similar to the case studied in \cite{Sottinen-Viitasaari-2017b} in the fractional Brownian setting. More precisely, we have
\begin{eqnarray*}
\hat m^{a,0}_t(u) &=& \frac{1}{a}\int_0^u k(t,s)\d W^{a,0}_s \\
&=& \frac{1}{a}\int_0^u k(t,s)\d (a W_s) \\
&=& \int_0^u k(t,s)\, \d W_s
\end{eqnarray*}
and
\begin{eqnarray*}
\hat r^{a,0}(t,s|u) &=&
\int_0^{t\wedge s} \left(1-\1_v(u)\right)^2k(t,v)k(s,v)\, \d v \\
&=&
\int_u^{t\wedge s} k(t,v)k(s,v)\, \d v \\
&=&
\int_0^{t \wedge s} k(t,v)k(s,v)\, \d v - \int_0^u k(t,v)k(s,v)\, \d v \\
&=&
r(t,s) - \int_0^u k(t,v)k(s,v)\, \d v.
\end{eqnarray*}
In particular, we see that $X_u$ is $\F^{W}_u$-observable, since $\hat r^{a,0}(u,u|u) = 0$.

The reason to use $a\ne1$ parametrization is that the general $(a,b)$ parametrization is easily translated into the correlation parametrization. Indeed, let $\rho\in[-1,1]$.  Then, setting $a=\rho$ and $b=\sqrt{1-\rho^2}$ we have that $W^{a,b}$ is a Brownian motion that is correlated with the driving Brownian motion $W$ with correlation coefficient $\rho$. 
\end{rem}

\begin{rem}[Predicting the Present]\label{rem:present}
Let us note that $X_u$ is not observable by the information $\F_u^{a,b}$, if $b>0$. Indeed, for simplicity and without any actual loss of generality, let us take $a=1$. Then Theorem \ref{thm:cond_law} states that
\begin{eqnarray*}
\hat m^{1,b}_u(u)
&=&
\frac{1}{1+b^2}\int_0^u k(u,s)\, \d W^{1,b}_s \\
&=&
\frac{1}{1+b^2}\int_0^u k(u,s)\, \d W_s + \frac{b}{1+b^2}\int_0^u k(u,s)\, \d \tilde W_s \\
&=&
\frac{1}{1+b^2} X_u + \frac{b}{1+b^2}\int_0^u k(u,s)\, \d\tilde W_s.
\end{eqnarray*}
Setting
$
\tilde{X}_u = \int_0^u k(u,s)\d\tilde W_s,
$
we note that $\tilde X$ is an independent copy of $X$, meaning that our best prediction for $X_u$ given $\F^{1,b}_u$ is the linear combination 
$$
\frac{1}{1+b^2}X_u + \frac{b}{1+b^2}\tilde{X}_u
$$
that is observable, while $X_u$ or $\tilde{X}_u$ are not. Similarly, we may compute the variance 
\begin{eqnarray*}
\hat v^{1,b}(u|u) &=& \hat r^{1,b}(u,u|u) \\
&=&
\int_0^u \left(\frac{b^2}{1+b^2}\right)^2 k(u,v)^2\, \d v + \frac{b^2}{(1+b^2)^2}\int_0^u k(u,v)^2\, \d v  \\
&=&
\frac{b^2(b^2+1)}{(1+b^2)^2}\int_0^u k(u,v)^2 \, \d v \\
&=&
\left(\frac{b}{1+b}\right)^2 \int_0^u k(u,v)^2 \, \d v.
\end{eqnarray*}	
In particular, we see that, as $b\to0$,
\begin{eqnarray*}
\hat m_u^{1,b}(u) &\to& X_u, \\
\hat v^{1,b}(u|u) &\to& 0.	
\end{eqnarray*}	
\end{rem}

\begin{proof}[Proof of Theorem \ref{thm:cond_law}]
First we note that the representation \eqref{eq:gvp} ensures that we can consider our Gaussian processes as random objects taking values in a separable Hilbert space.  This ensures the existence of regular conditional laws. Moreover, conditional Gaussian processes are Gaussian with random mean and deterministic covariance.  Indeed, see e.g. Janson \cite[Chapter 9]{Janson-1997}.  Consequently, to determine the regular conditional law of a Gaussian process, one only needs to calculate its conditional mean and conditional covariance.

Let us begin by calculating the conditional mean. To show that \eqref{eq:cond_mean} is indeed the conditional expectation of $X_t$ given $\F^{a,b}_u$, it is enough to show that
\begin{equation}\label{eq:mean_condition}
\E\left[\left(X_t - \int_0^u \hat k^{a,b}(t,s)\, \d W^{a,b}_s \right) W^{a,b}_v\right]
= 0
\end{equation}
for all $v\le u$.  Indeed, this means that the residual term
$$
\varepsilon_t^{a,b}(u) = X_t-\int_0^u \hat k^{a,b}(t,s)\, \d W^{a,b}_s
$$ 
is orthogonal to all the random variables $W^{a,b}_v$, $v\le u$.  By Gaussianity, this implies that $\varepsilon^{a,b}_t(u)$ is independent of the $\sigma$-algebra $\F^{a,b}_u$.  Since we have the decomposition
$$
X_t = \hat m_t^{a,b}(u) + \varepsilon^{a,b}_t(u),
$$ 
and $\hat m^{a,b}_t(u)$ is clearly $\F^{a,b}_u$-measurable, this implies that $\hat m^{a,b}_t(u) = \E[X_t\,|\, \F^{a,b}_u]$.

Now, the condition \eqref{eq:mean_condition} can be rewritten as
\begin{eqnarray*}
\lefteqn{\E\left[\int_0^t k(t,s)\, \d W_s\, (aW_v+b\tilde W_v)\right]} \\
&=&
\E\left[\int_0^u \hat k^{a,b}(t,s)\, \d(aW_s+b\tilde W_s)\, (aW_v+b\tilde W_v)\right].
\end{eqnarray*}
Since $W$ and $\tilde W$ are independent, this simplifies into
\begin{eqnarray*}
\lefteqn{a\E\left[\left(\int_0^t k(t,s)\, \d W_s\right) W_v\right]} \\
&=&
a^2\E\left[\left(\int_0^u \hat k^{a,b}(t,s)\,\d W_s\right) W_v\right]
+ b^2 \E\left[\left(\int_0^u \hat k^{a,b}(t,s)\,\d\tilde W_s\right) \tilde W_v\right],
\end{eqnarray*}
where we also used the bilinearity 
$$
\int_0^t f(s)\,\d(a W_s+b \tilde W_s) = a \int_0^t f(s)\, \d W_s + b \int_0^t f(s)\, \d\tilde W_s,
$$
for all $f\in L^2([0,t])$.

By using the It\^o isometry, we obtain the criterion
\begin{equation}\label{eq:cond_mean_criterion}
a\int_0^{t\wedge v} k(t,s)\, \d s
=
\left(a^2+b^2\right)\int_0^{u\wedge v} \hat k^{a,b}(t,s)\, \d s.
\end{equation}
We now have two cases: (i) $u\leq t$ (predicting the future), and (ii) $u>t$ (predicting the past).

Case (i): 
Since $v\le u$ and  $t\ge u$ we obtain by differentiating \eqref{eq:cond_mean_criterion} with respect to $v$ that
$$
a k(t,v) = \left(a^2 +b^2\right)\hat k^{a,b}(t,v)
$$
for Lebesgue almost everywhere. The formula \eqref{eq:cond_mean} follows from this, since if $k_1(t,s)=k_2(t,s)$ for almost every $s$, then for any Brownian motion $B$ we have 
$$
\int_0^t k_1(t,s)\, \d B_s = \int_0^t k_2(t,s)\, \d B_s
$$
almost surely.

Case (ii):
Now $u>t$ and $v\leq u$. Suppose then first that $v<t$. Then differentiating \eqref{eq:cond_mean_criterion} with respect to $v$ gives us (Lebesgue almost everywhere)
$$
\hat k^{a,b}(t,v) = \frac{a}{a^2+b^2}k(t,v)
$$
showing \eqref{eq:cond_mean} for $v< t$ precisely as in the case (i).

Let then $v\ge t$. Since
$$
(a^2+b^2)\int_0^v \hat k^{a,b}(t,s)\, \d s 
= 
(a^2+b^2)\int_0^t \hat k^{a,b}(t,s)\,\d s + (a^2+b^2)\int_t^v \hat k^{a,b}(t,s)\,\d s
$$
and
$$
(a^2+b^2)\int_0^t \hat k^{a,b}(t,s)\, \d s = a\int_0^t k(t,s)\, \d s,
$$
\eqref{eq:cond_mean_criterion} implies that
$$
\int_t^v \hat k^{a,b}(t,s)\, \d s = 0.
$$
Differentiating again with respect to $v$ gives us (Lebesgue almost everywhere)
$$
\hat k^{a,b}(t,v) = 0
$$
Finally, noting that  
$
k(t,v) = 0
$
for $v>t$, completes the proof of \eqref{eq:cond_mean}. 

In order to conclude the proof we need to prove \eqref{eq:cond_var}. But this is now straightforward. Indeed, the general theory tells us that
$$
\hat r^{a,b}(t,s|u) =
\E\left[\varepsilon_t^{a,b}(u)\varepsilon_s^{a,b}(u)\right], 
$$
where
\begin{eqnarray*}
\lefteqn{\varepsilon_t^{a,b}(u)}\\ &=&
X_t - \hat m_t^{a,b}(u) \\
&=&
\int_0^t k(t,v)\, \d W_v - \int_0^u \hat k^{a,b}(t,v)\, \d W^{a,b}_v \\
&=&
\int_0^t k(t,v)\, \d W_v - a\int_0^u \hat k^{a,b}(t,v)\, \d W_v
- b\int_0^u \hat k^{a,b}(t,v)\, \d\tilde W_v \\
&=&
\int_0^t k(t,v)\, \d W_v - \frac{a^2}{a^2+b^2}\int_0^u k(t,v)\, \d W_v
- \frac{ab}{a^2+b^2}\int_0^u k(t,v)\, \d\tilde W_v \\
&=&
\int_0^t \left(1 - \frac{a^2}{a^2+b^2}\1_u(v)\right)k(t,v)\, \d W_v
- \frac{ab}{a^2+b^2}\int_0^u k(t,v)\, \d\tilde W_v 
\end{eqnarray*}
Consequently, by the independence of $W$ and $\tilde W$, and the It\^o isometry, we see that
\begin{eqnarray*}
\lefteqn{\hat r^{a,b}(t,s|u)} \\
&=&
\E\left[\int_0^t \left(1 - \frac{a^2}{a^2+b^2}\1_u(v)\right)k(t,v)\, \d W_v
\int_0^s \left(1 - \frac{a^2}{a^2+b^2}\1_u(v)\right)k(s,v)\, \d W_v\right] \\
& & +
\E\left[\frac{ab}{a^2+b^2}\int_0^u k(t,v)\d\tilde W_v \frac{ab}{a^2+b^2}\int_0^u k(s,v)\, \d\tilde W_v\right] \\
&=&
\int_0^{t\wedge s} \left(1-\frac{a^2}{a^2+b^2}\1_u(v)\right)^2 k(t,v)k(s,v)\, \d v \\
& &
+
\frac{a^2b^2}{(a^2+b^2)^2}\int_0^u k(t,v)k(s,v)\, \d v,
\end{eqnarray*}
showing the validity of formula \eqref{eq:cond_var}.
\end{proof}

\section{Application to measurement errors}\label{sec:application}

Consider the noisy model
$$
X^b_t = \int_0^t k(t,s)\, \d W^{1,b}_s
= X_t + b\tilde X_t,
$$
where we have denoted 
$$
\tilde X_t = \int_0^t k(t,s)\, \d\tilde W_s.
$$
Our goal is to estimate $X$ from the noisy observations $X^b=X+b\tilde X$.  This means that we are predicting the present, see Remark \ref{rem:present}.

By using the observable $X^b$ as an estimator directly, we obtain the error
$$
\E\left[\left(X_t-X_t^b\right)^2\right]
=
b^2 \E[X_t^2],
$$
since $\tilde X$ is an independent copy of $X$.  Thus, if $|b|$ is big, using the observable $X^b$ may lead to huge estimation errors.
However, by using the conditional mean \eqref{eq:cond_mean} as an estimator, we have
\begin{eqnarray*}
\hat m_t^{1,b}(t) - X_t 
&=&
\frac{1}{1+b^2}\int_0^t k(t,s)\, \d W^{1,b}_s - X_t \\
&=&
\frac{1}{1+b^2}\left[X_t + b\tilde X_t\right] - X_t \\ 
&=&
\frac{b}{1+b^2}\tilde X_t -\frac{b^2}{1+b^2} X_t
\end{eqnarray*}
Since $X$ and $\tilde X$ are independent copies, we obtain the error
\begin{eqnarray*}
\E\left[{\left(\hat m_t^{1,b}(t)-X_t\right)}^2\right]
&=&
\E\left[{\left(\frac{b}{1+b^2}\tilde X_t -\frac{b^2}{1+b^2} X_t\right)}^2\right] \\
&=&
\frac{b^2}{(1+b^2)^2}\E[X_t^2] + \frac{b^4}{(1+b^2)^2}\E[X_t^2] \\
&=&
\frac{b^2(1+b^2)}{(1+b^2)^2}\E[X_t^2] \\
&=&
\frac{b^2}{1+b^2}\E[X_t^2] \\
&\le& \left(1 \wedge b^2\right)\E[X_t^2].
\end{eqnarray*}
Consequently, the variance of the error can be considerably reduced if one knows the Volterra kernel $k$ and the noise parameter $b$.  In particular, the variance of the error in the estimator $\hat m^{1,b}_t(t)$ is bounded in $b$ and always outperforms the variance of the error of the simple estimator $X^b_t$.

\bibliographystyle{siam}
\bibliography{../../pipliateekki}
\end{document}